\documentclass{amsart}
  \usepackage{amscd,amssymb,epsfig}
   \usepackage{epic,eepic}

                     \numberwithin{equation}{section}

                     \newtheorem{propo}{Proposition}[section]
                     \newtheorem{corol}[propo]{Corollary}
                     \newtheorem{theor}[propo]{Theorem}
                     \newtheorem{lemma}[propo]{Lemma}
                     \theoremstyle{definition}

                     \theoremstyle{remark}

             \newcommand{\CC}{\mathbb{C}}
             \newcommand{\QQ}{\mathbb{Q}}
                     \newcommand{\ZZ}{\mathbb{Z}}

                     \newcommand{\Hom}{\operatorname{Hom}}
             \newcommand{\Ker}{\operatorname{Ker}}

                     \newcommand{\id}{\operatorname{id}}

\def\1{\hbox{\rm\rlap {1}\hskip .03in{\rm I}}}

\def\id{{  \rm {id}}}

\def\orlin {G}

  \def\Hom{{ \rm  {Hom}}}

  \def\log {{ \rm  {log}}}

  \def\Aut {{ \rm  {Aut}}}

\def\mod {{ \rm  {mod}}}

  \def\dim  {{ \rm  {dim}}}

\def\mod {{ \rm  {mod}}\,}

\def\id {{ \rm  {id}}}

 \def\Ker {{ \rm  {Ker}}}
  \def\Im {{ \rm  {Im}}}

\def\Irr {{ \rm  {Irr}}}

\begin{document}
      \title{On certain enumeration problems in  two-dimensional topology}
                     \author[Vladimir Turaev]{Vladimir Turaev}
                     \address{%
              Department of Mathematics \newline
\indent  Indiana University \newline
                     \indent Bloomington IN47405 \newline
                     \indent USA}
                     \begin{abstract} We   announce a solution  to  several
 enumeration problems in topology of surfaces. This includes an
enumeration of  homotopy classes of sections of locally trivial
fiber bundles over surfaces and a computation of non-abelian
1-cohomology of surfaces.
                     \end{abstract}
                     \maketitle

   \section*{Introduction}\label{int}

We announce a solution  to  several enumeration problems in
two-dimensional topology. One problem
  deals  with an arbitrary   locally trivial fiber
bundle $p:E\to W$ over a closed connected oriented surface $W$. The
bundle $p$ may have sections, i.e.,  continuous mappings $s:W\to E$
such that $ps=\id_W$.
    If the   fiber
$F$ of~$p$ is path-connected and the group $\pi_1(F)$ is finite,
  then the  number of  sections of~$p$, considered
  up to homotopy and a natural action of $\pi_2(F)$, is finite. We give a  formula
  computing this number     in terms of certain
2-dimensional cohomology classes   associated with irreducible
complex linear representations of $\pi_1(F)$.   This yields, in
particular,   the following solution to the existence problem for
sections:  the bundle $p$ has a section if and only if the integer
number produced by our formula is non-zero. Note that the integer in
question is always non-negative and is computable provided one has
an efficient description of both the set of equivalence classes of
irreducible representations of $\pi_1(F)$ and the
  action of $\pi_1 (W)$ on this set determined by~$p$. As a specific
application, note the following theorem: in the case where the group
$\pi_1(F)$ is abelian (and finite), the bundle $p:E\to W$ has a
section if and only if the induced homomorphism $p_*:H_2(E)\to
H_2(W)$ is surjective. In particular, if  $p:E\to W$ is   a
principal $H$-bundle, where $H$ is a  connected topological group
with finite fundamental group, then $p$ is trivial if and only if
$p_*(H_2(E))= H_2(W)$.

The enumeration problem  for the sections of $ E\to W$   can be
 reformulated in terms of
  sections of the induced  homomorphism $ \pi_1(E)\to
\pi_1(W)$.   This generalizes to the following question. Given a
group epimorphism $ G'\to G$ with finite kernel~$\Gamma $ and a
homomorphism $g:\pi_1(W)\to G$, calculate the number of lifts of~$g$
to~$ G'$. The finiteness of  $\Gamma$ ensures that this number is
finite (possibly, zero). Our main result computes this number in
terms of   2-dimensional cohomology classes associated with
irreducible representations of~$\Gamma$. This encompasses the
problem of finding whether or not $g$ lifts to $ G'$.

 Other
enumeration problems considered here deal with   counting  principal
fiber bundles over $W$ and with a computation of non-abelian
1-cohomology of $W$.

Throughout this paper, we fix  two (discrete) groups $G$, $G'$  and
an epimorphism $q:G'\to G$ with finite kernel $\Gamma=\Ker\, q$. The
symbol $W$   denotes  a closed connected oriented surface of
positive genus with fundamental group $\pi$.

\section{Enumeration of homomorphisms}\label{eh}

  By a {\it
 (linear) representation}  of~$\Gamma$, we mean a   homomorphism
$ \Gamma \to  GL_n (\CC)$ with   $n=1,2,\ldots$.   Two
 representations $\rho :\Gamma\to GL_{n} (\CC)$ and
 $\rho':\Gamma\to GL_{n'}(\CC)$ are {\it equivalent}
  if $ n =n'$ and there
 is a  matrix $M\in GL_{n }(\CC)$ such that $\rho'(h)=M^{-1}\rho  (h) M $ for all $h\in
 \Gamma$. A
 representation $ \Gamma\to  GL_n (\CC)$ is {\it
irreducible}  if the only linear subspaces of $\CC^n$ preserved
under the induced action of $\Gamma$ are $0$ and $\CC^n$. The set of
equivalence classes of irreducible representations of $\Gamma$ is
denoted $\Irr (\Gamma)$.

Since $\Gamma=\Ker\,  q $ is a  normal subgroup of $G'$, the group
$G'$ acts on $\Gamma$ by conjugations. The induced  action of $G'$
on $\Irr (\Gamma)$  is trivial on $\Gamma\subset G'$ and therefore
induces an action of $G=G'/\Gamma$ on $\Irr (\Gamma)$. Given an
irreducible representation $\rho:\Gamma\to GL_n(\CC)$,
  let    $G_\rho\subset G$ be the stabilizer of the equivalence class of $\rho$
  under this action of $G$ on $\Irr (\Gamma)$. Thus,
 $G_\rho$  consists of  all $\alpha\in   G$ such that  for some   $\widetilde \alpha\in
 q^{-1} (\alpha)$, the   representation   $\Gamma \to GL_n(\CC)$,
 $ h\mapsto \rho (\widetilde \alpha^{-1} \, h \, \widetilde \alpha)$ is
 equivalent to~$\rho$.
Note that if the latter condition holds for some
 $\widetilde \alpha\in
 q^{-1} (\alpha)$, then it holds for all $\widetilde \alpha\in
 q^{-1} (\alpha)$. The subgroup    $G_\rho$ of $G$
depends only on the equivalence class of $\rho$.

The    representation $\rho $ determines a cohomology class
$\zeta_\rho \in H^2(G_\rho; \CC^*)$ as follows. For each $\alpha\in
G_\rho $, pick $\widetilde \alpha\in
 q^{-1} (\alpha)\subset G'$.
Then     there is a matrix $M_\alpha\in GL_n(\CC)$ such that
 $$\rho (\widetilde \alpha^{-1} \, h\, \widetilde
 \alpha )= M_\alpha^{-1}  \, \rho (h) \, M_\alpha    $$
 for all $h\in \Gamma$.
 The irreducibility of $\rho$ implies that
  $M_\alpha$ is
  unique  up to multiplication by an element of $\CC^*$; we fix $M_\alpha$ for all $\alpha\in G_\rho$.
For any $\alpha, \beta\in G_\rho$, we have
 $ \widetilde {\alpha\beta}^{-1}\, \widetilde \alpha \, \widetilde \beta
  \in \Gamma$. It is easy to check that
  there is a unique   $\zeta_{\alpha, \beta} \in
\CC^*$ such that \begin{equation}\label{1.22.a} \zeta_{\alpha,
\beta} \, M_\alpha\, M_\beta = M_{\alpha \beta}   \, \rho(\widetilde
{\alpha\beta}^{-1}\, \widetilde \alpha \, \widetilde \beta) \, .
\end{equation} The family $\{\zeta_{\alpha,\beta}\}_{\alpha,\beta}$
turns out to be a
 2-cocycle on
 $G_\rho$. Its cohomology class  $\zeta_\rho\in H^2(G_\rho; \CC^*)$
 depends only
  on the equivalence
class of $\rho$ and does not depend on the choice
 of   $\{M_\alpha\}_{\alpha}$  or   $\{\widetilde \alpha\}_{\alpha
 }$. (The definitions   of $G_\rho$ and
  $\zeta_\rho$    do   not use the finiteness of
$\Gamma$.)

Given a   homomorphism $g:\pi=\pi_1(W)\to  G$, a {\it lift}  of $g$
to $G'$ is a homomorphism $g': \pi \to G'$ such that $qg'=g$. The
set (possibly empty) of all such lifts is denoted by $\Hom_g (\pi
,G')$. Since $\pi $ is finitely generated and $\Gamma$ is finite,
the  set $\Hom_g (\pi ,G')$ is finite. The number of its elements is
bounded from above by $ {\vert \Gamma \vert}^b$, where $b$ is the
first Betti number of $W$ and the vertical bars stand for the
cardinality of a set. We now compute $\vert \Hom_g (\pi ,G')\vert$
 in terms of representations of
$\Gamma$. Note that the finiteness of $\Gamma$   guarantees  that
the set $\Irr (\Gamma)$ is finite.

 \begin{theor}\label{t1}    Let
$g:\pi=\pi_1(W)\to G$ be a  group homomorphism. Then
\begin{equation}\label{1.2.a} \vert  \Hom_g (\pi ,G')\vert\, =\, \vert  \Gamma \vert  \,
\sum_{\rho\in \Irr(\Gamma),\, G_\rho \supset \,g(\pi)  } \left
(\frac{\vert \Gamma \vert}{\dim \, \rho}\right )^{-\chi(W)}
\,g^*(\zeta_\rho) ([W]) ,\end{equation}  where $\rho$ ranges over
the equivalence classes of irreducible (complex)
    representations of $\Gamma$  such that
$G_\rho \supset g(\pi)$, and $g^*(\zeta_\rho) ([W])\in \CC^*$ is the
evaluation of  $g^*(\zeta_\rho) \in H^2(\pi; \CC^*)$ on the
fundamental class $[W]\in H_2(W{{}})=H_2(\pi{{}})$   of $W$.
\end{theor}

Here and below the unspecified group of coefficients in homology
is $\ZZ$. The cohomology class $g^*(\zeta_\rho) $ is
 well-defined because    $g(\pi)\subset
G_\rho$.  The addition on the right-hand side of (\ref{1.2.a}) is
the addition in $\CC$. The sum on the right-hand side of
(\ref{1.2.a}) is always non-empty because the trivial
one-dimensional representation $\rho_0:\Gamma\to \{1\}\in \CC^*=
GL_1(\CC)$ satisfies $G_{\rho_0}=G\supset g(\pi)$ and
$g^*(\zeta_{\rho_0}) ([W])=1$.

  Note a few cases where Formula (\ref{1.2.a}) is known or may be directly deduced from the known results.       If $g=1$, then
 $\Hom_g (\pi ,G')=\Hom (\pi ,\Gamma) $. In this case $G_\rho\supset g(\pi)$ and $g^*(\zeta_\rho)([W])=1$ for
all irreducible  representations $\rho$ of $\Gamma$. Formula
(\ref{1.2.a}) is then equivalent to the well-known Frobenius-Mednykh
formula
\begin{equation}\label{1.1.a}\vert  \Hom (\pi , \Gamma)\vert\, =\,
\vert  \Gamma \vert  \, \sum_{\rho\in \Irr(\Gamma)} \,\left
(\frac{\vert \Gamma \vert}{\dim \, \rho}\right )^{-\chi(W)}\,
.\end{equation}
  For
$W=S^2$, this formula   is equivalent to the classical equality
$\sum_\rho (\dim \, \rho)^2=\vert \Gamma \vert$. For $W=S^1\times
S^1$, Formula (\ref{1.1.a}) was first established by Frobenius
\cite{Fr}. The general case of (\ref{1.1.a}) is due to Mednykh
\cite{Me}, see also \cite{FQ} and \cite{Jo}.

If $G'=\Gamma\times G$ and $q:G'\to G$ is the projection, then
(\ref{1.2.a})   directly follows from (\ref{1.1.a}) since in this
case $\Hom_g (\pi ,G')=\Hom  (\pi ,\Gamma)$ and $\zeta_\rho$ is
trivial for all $\rho$.

If
  $\Gamma=\Ker\, q$ is central in $G'$, then Formula (\ref{1.2.a})  is
essentially obvious. In  this case $G_\rho=G$ and $\dim \, \rho=1$
for all~$\rho$,  while $\zeta_\rho$ is the image of the standard
cohomology class $\zeta\in H^2(G; \Gamma)$ determined by $q$ under
the coefficient homomorphism $H^2(G; \Gamma)\to H^2(G; \CC^*)$
induced by $\rho:\Gamma\to \CC^*$. Formula (\ref{1.2.a}) can be
deduced then from the following   easy assertions: $g$ has a lift to
$G'$ if and only if $g^*(\zeta) =1$; if there are such lifts, then
their number is equal to $ \vert \Gamma \vert^{ 2-\chi(W)} $.

A proof   of Theorem \ref{t1} in the full generality uses the
techniques of topological quantum field theory and  will be given
elsewhere. Note that Formula~(\ref{1.1.a}) extends to surfaces with
boundary and to non-orientable surfaces, see \cite{Jo}, \cite{Sn}.
Theorem~\ref{t1} admits similar extensions, but we shall not discuss
them here.

 We keep the assumptions of
 Theorem \ref{t1} and establish several corollaries. Observe that $$g^*(\zeta_\rho) ([W])= (\zeta_\rho
\vert_{g(\pi)}) (g_*([W]))\, , $$ where  $\zeta_\rho \vert_{g(\pi)}
\in H^2(g(\pi);\CC^*)$ is the restriction of $\zeta_\rho\in
H^2(G_\rho;\CC^*)$ to $g(\pi) $ and   $$g_*:
H_2(W{{}})=H_2(\pi{{}})\to H_2 (g(\pi){{}})$$ is the homomorphism
induced by $g$. Formula (\ref{1.2.a}) can   be rewritten as
\begin{equation}\label{1.3.a} \vert \Hom_g (\pi ,G')\vert
  =\,  \vert  \Gamma \vert \, \sum_{\rho\in \Irr(\Gamma),\, G_\rho \supset
\,g(\pi)  } ( \vert  \Gamma \vert/\dim \, \rho)^{
-\chi(W)}\,(\zeta_\rho \vert_{g(\pi)}) (g_*([W]))\, .\end{equation}
This implies the following claim.

    \begin{corol}\label{1.3.1}    The number $\vert  \Hom_g (\pi ,G')\vert $ is
     determined by the homomorphism $q:G'\to G$,  the
genus of $W$, the  group $  g(\pi)\subset G$,  and the homology
class $ g_*([W])\in H_2 (g(\pi){{}})  $. \end{corol}

Note the following   special case of (\ref{1.3.a}).

\begin{corol}\label{1.3.1+} If $g:\pi\to G$ is an epimorphism, then
\begin{equation}\label{1.3.a+}\vert  \Hom_g (\pi ,G')\vert\,=\, \vert  \Gamma \vert  \,
\sum_{\rho\in \Irr(\Gamma),\, G_\rho =G  } (\vert  \Gamma \vert/\dim
\, \rho)^{-\chi(W)}\, \zeta_\rho(g_*([W]))\, .\end{equation}
\end{corol}

To proceed, we need the following property of   $\zeta_\rho$.
Consider the group $G'_\rho=q^{-1}(G_\rho)\subset G'$. Taking
determinant  on both sides  of (\ref{1.22.a}),  one easily obtains
that \begin{equation}\label{1.3.ddd} q^*((\dim\, \rho) \,
\zeta_\rho)=0\, , \end{equation} where $q^*:H^2(G_\rho;\CC^*)\to
H^2( G'_\rho;\CC^*)$ is the homomorphism induced by~$q$. This
implies that the values of $\zeta_\rho$ on the image of the
homomorphism $q_*:H_2(G'_\rho )\to H_2( G_\rho )$ are roots of unity
of order $\dim \, \rho$. The finiteness of $\Gamma$ implies that the
image of $q_*$ is a subgroup of finite order in $H_2( G_\rho )$.
Hence the values of $\zeta_\rho$ on all elements of $ H_2( G_\rho )$
are roots of unity and, in particular, have absolute value   1. Now,
comparing (\ref{1.3.a}) with (\ref{1.1.a}) termwise, we obtain the
following inequality.

    \begin{corol}\label{1.3.3}  For any homomorphism
    $g:\pi=\pi_1(W)\to G$,
    \begin{equation}\label{1.3.d}
    \vert  \Hom_g (\pi ,G')\vert \leq \vert  \Hom  (\pi ,\Gamma)\vert\, .\end{equation} This
inequality is an equality if and only if $g(\pi)\subset G_\rho$ and
$g^*(\zeta_\rho) ([W])=1$ for all irreducible   representations
$\rho$ of $\Gamma$.\end{corol}

The inequality (\ref{1.3.d}) does not hold for arbitrary groups
$\pi$. For example, let  $G'=S_3$, the group of permutations of the
set $\{1,2,3\}$,   and let $q:G'\to {\ZZ}/2\ZZ$ be the
  surjection sending all transpositions to $1\, (\mod \, 2)$.
  Clearly, $\Gamma= \Ker\, q={\ZZ}/3\ZZ$.
Let $\pi$ be the group with $m\geq 1$ generators $x_1$, $\ldots$,
$x_m$ and defining  relations $x^2_1=x^2_2=\cdots =x^2_m$. A
homomorphism $\pi\to \Gamma$ has to  send all the generators $x_1$,
$\ldots$, $x_m$ to the same element. Therefore, $\vert  \Hom_g (\pi
,\Gamma)\vert=3$. On the other hand, the homomorphism $ \pi\to
{\ZZ}/2\ZZ$ sending $x_1$, $\ldots$, $x_m$ to $1\, (\mod \, 2)$
admits at least $3^m$ lifts to $G'$ sending $x_1$, $\ldots$, $x_m$
to arbitrary transpositions.

A {\it section}\index{homomorphism!section of}  of a group
homomorphism $p:\pi' \to \pi$ is a homomorphism $s:\pi\to \pi'$ such
that $ps=\id_\pi$. The set of sections of $p$ is denoted
by~$S_*(p)$.

   \begin{corol}\label{1.3.4}    Let $p:\pi' \to \pi=\pi_1(W)$ be a group epimorphism with finite kernel $\Phi$.
    Then
 $$\vert S_*(p) \vert=\vert  \Phi \vert \,
\sum_{\rho\in \Irr(\Phi),\, \pi_\rho = \pi   } (\vert  \Phi
\vert/\dim \, \rho)^{-\chi(W)}\, \zeta_\rho  ([W]) \, . $$
 \end{corol}

    This    is obtained from  Theorem \ref{t1} by setting
$G=\pi$, $G'=\pi'$, $\Gamma=\Phi$, $q=p$, and $g=\id:\pi\to \pi$.
   Corollary \ref{1.3.3} implies that
 $ \vert S_*(p) \vert\leq \vert \Hom (\pi ,\Phi)\vert $.  This
inequality is an equality if and only if $\pi_\rho=\pi$ and $
\zeta_\rho ([W])=1$ for all irreducible
  representations $\rho$ of $\Phi$. The results of the next section  imply that
   the number $ \vert S_*(p) \vert$ is divisible by $\vert \Phi \vert \,
\vert Z(\Phi)\vert^{2d-2}$, where~$Z(\Phi)$
 is the center of~$\Phi$ and $d$ is the genus of $W$.

\section {The functions $\{v_k\}_k$}\label{vvv}

 The aim of this section is to deduce from Theorem \ref{t1} the following claim.

 \begin{theor}\label{1.3.2}   Let
$g:\pi=\pi_1(W)\to G$ be a  group homomorphism.  Then the number
$\vert \Hom_g (\pi ,G')\vert$ is divisible by $\vert \Gamma \vert \,
\vert Z(\Gamma)\vert^{2d-2}$, where~$Z(\Gamma)$
 is the center of~$\Gamma$ and $d$ is the genus of $W$.\end{theor}

Observe first that for each $k=1,2, \ldots$, the
   epimorphism $q:G'\to G$  determines a function
   $v_k:H_2(G{{}})\to \CC$ by
$$v_k (h)=\sum_{\rho\in \Irr(\Gamma),\,  \dim\, \rho=k, \,G_\rho=G} \,\,\zeta_\rho (h)\, \in \CC\,  ,$$
where  $h\in H_2(G{{}})$.
 It is well known that the dimension of any
irreducible representation    of $\Gamma$ divides $\vert
\Gamma/Z(\Gamma) \vert$ and is smaller than or equal to
  $ {\vert
\Gamma/Z(\Gamma)\vert}^{1/2}$. Therefore $v_k=0$ if $k$ does not
divide $\vert \Gamma/Z(\Gamma) \vert$ or
  $k>  {\vert
\Gamma/Z(\Gamma)\vert}^{1/2}$.

We can  rewrite Formula (\ref{1.3.a+})   in terms of the   functions
$v_1, v_2,\ldots$: for any epimorphism $g:\pi=\pi_1(W)\to G$,
\begin{equation}\label{vv1}\vert \Hom_g (\pi ,G')\vert =\, \vert
\Gamma \vert
 \, \sum_{k\geq 1  }   \,v_k (g_*([W]))\, ({\vert
\Gamma \vert/k})^{-\chi(W)}\,
  . \end{equation}

The next  lemma summarizes the properties of the functions
$\{v_k\}_k$.

  \begin{lemma}\label{l5}
   For $k=1,2,\ldots$, let $N_k$ be  the number of  equivalence classes of  irreducible $k$-dimensional
  complex representations $\rho$ of $\Gamma$ such that $G_\rho=G$.
  Let $Q  $ be the image of the homomorphism $
q_*: H_2(G'{{}})\to H_2(G{{}})$. Then

(a) For all $k$, the function   $v_k$ takes  only integer values and
is zero outside~$Q$. Moreover, all values of $v_k$ lie in the set
$\{-N_k, -N_k+1, \ldots, N_k-1, N_k\}$;

(b) For all
  $h\in Q$, we have  $v_1(h)=N_1=\vert
\Gamma/[\Gamma, G'] \vert$;

(c)  For all $k $, we have $v_k(0)=N_k$.
\end{lemma}

    \begin{proof} (a)  It is well known that for any   $h\in H_2(G{{}})$ there are a
 closed connected oriented  surface $\Sigma$ and  a
  homomorphism $g:\pi_1(\Sigma)\to  G$ such that
  $g_*([\Sigma])=h$. We   say that the pair $(\Sigma,g)$ {\it realizes} $h$.
  Pick a representation  $\rho$ of $\Gamma$ such that $  G_\rho \neq G$.
    Adding to $\Sigma$ a handle   and mapping its meridian  to $1\in G$
  and its longitude  to any element  of $G-G_\rho$,
  we obtain a realization of $h$ by a    surface $\Sigma'$ and  a
  homomorphism
  $\pi_1(\Sigma')\to G$  whose image
 meets $G-G_\rho$. Repeating this process, we can
   realize $h$ by a   surface $\Sigma_+$ with fundamental group $\pi$ and  a
  homomorphism $g:\pi \to G$  such that $g(\pi)\subset G_\rho$ only when $G_\rho=G$. Formula (\ref{1.3.a}) implies then that
\begin{equation}\label{vv2}\vert  \Hom_g (\pi ,G')\vert     =\,  \vert  \Gamma \vert^{2d-1}
 \, \sum_{k\geq 1  }   \,v_k (h)\,  {k}^{2-2d}
  , \end{equation}
  where $d$ is the genus of $\Sigma_+$.
 Any surface   of   bigger genus    admits a degree
  one map to~$\Sigma_+$. Such a map induces a surjection of fundamental
  groups. We can apply (\ref{vv2}) to the composition of this surjection
  with $g$. This
  implies    that   $\sum_{k\geq 1  }  \,v_k (h)\,k^{2-2n}\in \QQ$  for all   $n\geq d$.
  By   linear
  algebra, $v_k (h)\in \QQ$ for all $k$.

By the remarks made before the statement of Corollary \ref{1.3.3},
for any irreducible
   representation  $\rho$ of $\Gamma$ with $G_\rho=G$ and any $h\in H_2(G{{}})$,  the   number
  $\zeta_\rho (h)\, \in \CC $
  is a  root  of unity. Thus,  $v_k(h)$ is a sum of
 $N_k$ roots of unity. Hence,  $v_k(h)$ is
     an algebraic integer. Therefore $ v_k(h)\in \ZZ$ and $\vert v_k(h)
\vert \leq N_k$.

  If $h\notin Q$, then a   homomorphism  $g $
  realizing $h$ as above cannot lift  to~$G'$. Formula  (\ref{vv2}) and the
  argument after this formula
show  that $\sum_{k\geq 1  }  \,v_k (h)\,k^{2-2 n }=0$ for all
sufficiently big natural numbers~$n$. This gives a non-degenerate
system of linear equations on  $ \{ v_k(h)\}_k $.   Therefore
$v_k(h)=0$ for all~$k$.

(b) By definition, $v_1(h)=\sum_\rho \zeta_\rho (h)$, where $\rho$
runs over all homomorphisms $ \Gamma \to \CC^*$ such that
$G_\rho=G$. Formula  (\ref{1.3.ddd})  and the inclusion $h\in Q$
imply that $\zeta_\rho (h)=1$. Therefore $v_1(h)=N_1$ is simply the
number of such $\rho$.
   The condition $G_\rho=G$ holds  if and only if
   $\rho(aha^{-1}h^{-1})=1$ for all $a\in G'$, $h\in \Gamma$. The latter holds if and only if
   $\rho ([\Gamma, G'])=1$. Thus,
   $N_1=\vert \Gamma / [\Gamma, G']\vert $.

(c) The equality $v_k(0)=N_k$ follows from the equality
  $\zeta_\rho (0)=1$ for  all irreducible
  representations $\rho$ of $\Gamma$.
  \end{proof}

As an exercise, the reader may verify that   $v_k(-h)= {v_k(h)}$ and
$v_k(h+kh')=v_k(h)$ for all $h\in H_2(G{{}})$ and  $h'\in Q$.

\subsection* {Proof of Theorem \ref{1.3.2}.} Replacing
  $G$ and  $G'$  by   $g(\pi)$ and   $q^{-1}(g(\pi))$, respectively,   we
can reduce ourselves to the case where $g$ is an epimorphism. In
this case Theorem \ref{1.3.2} follows   from Formula (\ref{vv1}),
Lemma \ref{l5}(a), and the fact that $v_k=0$   if $k$ does not
divide $\vert \Gamma/Z(\Gamma)\vert=\vert \Gamma\vert/ \vert
Z(\Gamma)\vert$.

\section {The homological  obstruction to lifting}

  Consider in more detail
  the question of the existence of   lifts to $G'$ for a given homomorphism
$g:\pi=\pi_1(W) \to G$. Replacing
  $G$ and  $G'$  by   $g(\pi)$ and   $q^{-1}(g(\pi))$, respectively,   we
can reduce ourselves to the case where $g$ is an epimorphism. Then
Formula (\ref{vv1}) computes $\vert \Hom_g (\pi ,G')\vert$ in terms
of the numbers $\{v_k(h)\}_k$, where $h=g_*([W])\in H_2(G)$. Thus,
$g$ lifts to $G'$ if and only if the right-hand side of (\ref{vv1})
is non-zero.

Another approach to the same question stems from   homological
considerations. If an epimorphism $g:\pi \to G$ lifts to $G'$, then
the homology class $g_*([\Sigma]) \in H_2(G{{}})$ necessarily lies
in the image of the homomorphism $q_*: H_2(G'{{}})\to H_2(G{{}})$.
If the latter condition is satisfied, then we   say that the
homological obstruction to the lifting of $g$ to~$G'$ vanishes. In
general, the vanishing of the homological obstruction may not imply
that $g$ lifts to $G'$. The next two theorems show that there are no
further obstructions if the group $\Gamma=\Ker \, q$ is abelian or
the genus of $W$ is big enough.

In the sequel, the symbol
  $[\Gamma,
 G']$ denotes the subgroup of $\Gamma$ generated by the commutators of elements of $\Gamma$
  with elements of $G'$.

   \begin{theor}\label{7.1.1}
      Suppose that the group  $\Gamma $ is abelian. An epimorphism $g:\pi=\pi_1(W)\to G$  lifts
to $G'$ if and only if the homological  obstruction to the lifting
  vanishes. Moreover, if $g$ lifts
to $G'$, then \begin{equation}\label{eqq1}\vert \Hom_g (\pi ,
G')\vert = \vert \Gamma \vert^b\, \vert   [\Gamma,
 G']\vert^{-1}\, ,\end{equation}
  where $b=2-\chi(W)$ is the first Betti number of $W$.
    \end{theor}

\begin{proof}  Assume that   the
homological obstruction in question vanishes so that $h=g_*([W])\in
\Im\, q_*$. By Lemma \ref{l5}(b), we have $v_1(h)=\vert
\Gamma/[\Gamma, G'] \vert $. Since $\Gamma $ is abelian, all
irreducible   representations of $\Gamma$ are one-dimensional. So,
$v_k(h)=0$ for   $k\geq 2$.   Now, the claim of the theorem directly
follows from (\ref{vv1}).\end{proof}

      Formula (\ref{eqq1}) can be rewritten as
\begin{equation}\label{eqq2} \vert   \Hom_g (\pi , G')\vert=
    \vert \Hom  (\pi , \Gamma)\vert  \times  \vert   [\Gamma, G']\vert^{-1}\, . \end{equation}
This formula
  does not directly extend to   groups $\pi$ distinct from the fundamental groups of closed oriented surfaces.
  For instance, if $\pi$ is a free group of rank $n$, then
  the left-hand and right-hand sides of (\ref{eqq2}) are equal respectively
  to $\vert \Gamma \vert^n$ and
  $\vert \Gamma \vert^n\, \vert   [\Gamma, G']\vert^{-1}$. These   numbers are equal if and only if
  $[\Gamma, G']=1$, i.e.,  if and only if $\Gamma$   lies in the center of $G'$.

 \begin{theor}\label{7.1.2}      Suppose that    the first Betti number $b=2- \chi (W)$  of~$W$
satisfies
\begin{equation}\label{eqq3}b>\log_2(\vert [\Gamma, G']\vert -1)\, .\end{equation}
 An epimorphism
$g:\pi=\pi_1(W)\to G$  lifts to $G'$ if and only if the homological
obstruction to the lifting   vanishes. Moreover, if $g$ lifts to
$G'$, then
 \begin{equation}\label{eqq4}\vert \Hom_g (\pi , G')\vert\geq  \vert \Gamma \vert^b\, \vert   [\Gamma,
 G']\vert^{-1} \times (1-\frac{\vert [\Gamma,
  G']\vert-1}{2^b})
 \, .\end{equation}
     \end{theor}

   It is understood that if $[\Gamma, G']=1$, then the condition   (\ref{eqq3}) is empty.

\begin{proof} For $k=1,2,
\ldots$,
 denote by $M_k$ the number of  equivalence classes of  irreducible
$k$-dimensional
  representations   of $\Gamma$.
  It is clear that
  $\sum_{k\geq 1} \, M_k k^2=\vert \Gamma\vert$.

  Pick any $h\in H_2(G{{}})$. Let $N_k$ be  the same number as in Lemma \ref{l5}. The inequality $\vert v_k (h)\vert   \leq N_k  $
   established in Lemma \ref{l5}  and the obvious inequality $N_k\leq M_k$ imply
   that for any integer $n<2$,
  $$   \sum_{k\geq 1  }  \,v_k (h)\,k^{  n}
  \geq v_1(h) - \sum_{k\geq 2  }  \,\vert v_k (h)\vert \,k^{
   n}
  \geq v_1(h) - \sum_{k\geq 2  }  \,M_k \,k^{  n}
 $$
  $$= v_1(h) - \sum_{k\geq 2  }  \,M_k k^2 \,k^{n-2  }\,
 \geq v_1(h) - (\sum_{k\geq 2  }  \,M_k k^2) 2^{n-2  }
 $$
  $$= v_1(h) - 2^{n-2  }\, (\vert \Gamma\vert -M_1) \geq v_1(h) - 2^{n-2  }\, (\vert \Gamma\vert -N_1)\, .$$
 If $h\in  \Im\, q_*$, then   $v_1(h)=N_1=\vert \Gamma  /  [\Gamma, G']\vert$
  by Lemma \ref{l5}.
 This gives
$$  \sum_{k\geq 1  }  \,v_k (h)\,k^{  n}
 \geq \vert \Gamma/[\Gamma, G']\vert \left (1- \frac { \vert  [\Gamma, G']\vert
 -1 }{2^{2-n}} \right )\, .$$
 Setting   $n=\chi(W)=2-b$ and combining   with (\ref{vv1}), we  obtain (\ref{eqq4}).
 \end{proof}

\section {Geometric applications of Theorem  \ref{t1}}\label{kol}

\subsection* {Enumeration of
principal fiber bundles.}  Let ${\mathcal P}={\mathcal P}(W,
\Gamma)$ be the set of isomorphism classes of principal
  $\Gamma$-bundles over $W$. Recall that a
homomorphism
  $g:\pi
  \to
  \Gamma$ determines a principal
  $\Gamma$-bundle  $\xi_g$ over $W$, and every principal
  $\Gamma$-bundle   over $W$ is isomorphic to $\xi_g$ for some    $g $.
    Two   homomorphisms $g_1, g_2:\pi
  \to
  \Gamma$ determine isomorphic principal
  $\Gamma$-bundles over $W$ if and only if $g_1=hg_2 h^{-1}$ for some $h\in
  \Gamma$. Therefore
   ${\mathcal P}=  \Hom (\pi  , \Gamma)/\Gamma $,
  where $\Gamma$ acts on $\Hom (\pi  , \Gamma)$ by  conjugation. The stabilizer of $g\in \Hom (\pi  , \Gamma)$
  under this action is the   group $\{h\in \Gamma\, \vert \,
   hgh^{-1}=g\}$  isomorphic to   the group
  of automorphisms $\Aut (\xi_g )$ of $\xi_g$.
Combining these facts, we obtain
  $$\vert  \Hom (\pi  , \Gamma)\vert \, = \, \sum_{\xi \in {\mathcal P} } \,\frac{\vert
  \Gamma\vert }{\vert \Aut (\xi )\vert} \, . $$
   The Frobenius-Mednykh formula (\ref{1.1.a}) implies therefore that
\begin{equation}\label{ffm}\sum_{\xi\in {\mathcal P}(W; \, \Gamma)} \,\frac{1}{\vert \Aut
(\xi )\vert} \, = \,
  \sum_{\rho\in \Irr(\Gamma)}  \,
\left ( \frac {\vert  \Gamma \vert} {\dim \, \rho} \right
)^{-\chi(W)}\,  .\end{equation}

Theorem \ref{t1} yields a relative version of (\ref{ffm}) as
follows. Fix a principal $G$-bundle $\xi$ over $W$. By a {\it lift
of  $\xi$ to a principal $G'$-bundle}, we mean a pair (a principal
$G'$-bundle $\xi'$ over $W$,  an isomorphism of $G$-bundles
$f:\xi'/\Gamma \cong \xi$). Here $\xi'/\Gamma$ is the principal
$G$-bundle over $W$ obtained by factorizing the total space of
$\xi'$ by $\Gamma$. An {\it isomorphism} $(\xi'_1, f_1)\approx
(\xi'_2, f_2)$ of two such lifts of $\xi$ is an isomorphism of
principal $G'$-bundles $\xi'_1\to \xi'_2$ such that the induced
isomorphism of principal $G$-bundles $\xi'_1/\Gamma\to
\xi'_2/\Gamma$ composed with $f_2: \xi'_2/\Gamma \to \xi$ gives
$f_1$. In particular, an {\it automorphism} of a lift $ (\xi', f )$
of $\xi$ is an automorphism of  $\xi' $ inducing the identity on
$\xi' /\Gamma $. Such automorphisms form a group denoted by
$\Aut_\xi (\xi')$. The set of isomorphism classes of lifts of $\xi$
to principal $G'$-bundles is denoted by ${\mathcal P}(\xi)$.

Fix a homomorphism $g:\pi\to G$ such that $\xi = \xi_g$. It is clear
that every homomorphism
  $g'\in \Hom_g(\pi, G')$ determines a lift
of  $\xi$ to a principal $G'$-bundle, and every such lift arises
from some $g'\in \Hom_g(\pi, G')$.
    Two   homomorphisms $g'_1, g'_2\in \Hom_g(\pi, G')$ determine isomorphic lifts
of  $\xi$ if and only if $g'_1=hg'_2 h^{-1}$ for some $h\in
  \Gamma$. Therefore
  $${\mathcal P}(\xi)=  \Hom_g (\pi  , G')/\Gamma\, ,$$
  where $\Gamma$ acts on $\Hom_g (\pi  , G')$ by  conjugation. The stabilizer of any
  homomorphism  $g'\in \Hom_g (\pi  , G')$
  under this action is the   group $\{h\in \Gamma\, \vert \,
   hg'h^{-1}=g'\}$  isomorphic to    $\Aut_\xi (\xi_{g'} )$.
Combining these facts, we obtain
\begin{equation}\label{ffmb}\vert  \Hom_g (\pi  , G') \vert = \, \sum_{\xi' \in {\mathcal P}(\xi)} \,\frac{\vert
  \Gamma\vert }{\vert \Aut_\xi (\xi' )\vert} \, . \end{equation}
   Theorem \ref{t1} implies   that
\begin{equation}\label{ffmc}\sum_{\xi'\in {\mathcal P}(\xi)} \,\frac{1}{\vert \Aut_\xi (\xi' ) \vert} \,
  =\,
\sum_{\rho\in \Irr(\Gamma),\, G_\rho \supset \,g(\pi)  } \left (
\frac {\vert \Gamma \vert} {\dim \, \rho}\right
)^{-\chi(W)}\,\,g^*(\zeta_\rho) ([W]) \, .
\end{equation}

  Formula (\ref{ffmb}) and   Theorem \ref{1.3.2} imply that
    $\sum_{\xi'\in {\mathcal P}(\xi)} \, {1}/{\vert \Aut_\xi (\xi' ) \vert}$ is an integer
     divisible by $\vert Z(\Gamma)\vert^{2d-2}$,
   where   $d$ is the genus of $W$. Corollary
  \ref{1.3.3} gives
   $$\sum_{\xi'\in {\mathcal P}(\xi)} \, {1}/{\vert \Aut_\xi (\xi' ) \vert}\,\, \leq \,\,\vert \Gamma\vert^{-1}\,
   \vert \Hom (\pi, \Gamma)\vert\, .$$

Formula (\ref{ffmc})  implies that $\xi$ lifts to a principal
$G'$-bundle  if and only if the right hand side of (\ref{ffmc}) is
non-zero. Theorems \ref{7.1.1} and \ref{7.1.2} show that in the case
where   $\Gamma$ is abelian or the genus of $W$ is big enough, $\xi$
lifts to a principal $G'$-bundle if and only if the homology class
$g_*([W])\in H_2(g(\pi))$ lies in the image of the homomorphism
$q_*:H_2(q^{-1}(g(\pi)))\to H_2(g(\pi))$.

\subsection* {Enumeration of sections.} Theorem \ref{t1} may   be
used to count homotopy classes of sections of locally trivial fiber
bundles over the surface~$W$.     Let $p:E\to W$ be a locally
trivial fiber bundle with  fiber $F$.
  A {\it section}\index{section} of~$p$ is
a continuous mapping $s:W\to E$ such that $ps=\id_W$.   Two sections
  of $p$ are {\it homotopic} if they can be
  deformed into each other in the class of sections
  of $p$.   We say that two sections $ W\to E$ are obtained from each
other  by {\it bubbling}\index{bubbling} if they   coincide on the
complement of a small open disc $D\subset W$. The restrictions of
such two sections on the closed disc $\overline D\subset W$   form
then a mapping $S^2\to E$, a \lq\lq bubble". Two sections of $ p$
are
     {\it
 bubble equivalent}\index{bubble equivalence}   if they may be obtained from each other by
a finite sequence of   bubblings.  Decomposing a  deformation of a
section  into   local deformations, one easily observes that
homotopic sections are bubble equivalent. (If $\pi_2(F)=0$, then the
converse is also true so that the bubble equivalence is just the
homotopy.) Denote the set of bubble equivalence classes of  sections
of $p$ by ${\mathcal S}(p)$. We shall count the elements of this set
with certain weights, see Formula (\ref{9.2.c}) below.

The definition of  ${\mathcal S}(p)$ has a pointed version  as
follows. Fix a
 base point $e\in E$ and set $w=p(e)\in W$.
 A section
$s:W\to E$ of $ p$  is {\it pointed}\index{section!pointed} if
$s(w)=e$. Two  pointed sections  of $p$ are {\it homotopic} if they
can be
  deformed into each other in the class of   pointed
sections  of $p$. The definitions of the bubbling and   the bubbling
equivalence extend to   pointed sections in the obvious way with the
only difference that the disk $D$  in the definition of a bubbling
should lie in $W-\{w\}$. As above, homotopic pointed sections are
bubble equivalent, and the converse is true if $\pi_2(F)=0$. We
denote the set of bubble equivalence classes of pointed sections of
$p$ by ${\mathcal S}_* (p)$.

Suppose from now on that the fiber $F=p^{-1}(w)$ of $p$ is
path-connected. Set $\pi=\pi_1(W,w)$, $\pi'=\pi_1(E,e)$, and $\Phi
=\pi_1(F,e)$. The exact homotopy sequence of $p$ shows that the
homomorphism  $p_\#: \pi'\to \pi  $  is surjective and $\Ker\,
p_\#=\Phi$. By Section  1, every irreducible
   representation  $\rho$ of~$\Phi$ determines a subgroup $\pi_\rho$ of $\pi$ and a cohomology class $\zeta_\rho \in H^2(\pi_\rho;\CC^*)$.

  For any pointed
section $s:W\to E$ of $p$,   the induced homomorphism $s_\#: \pi \to
\pi'$
  is a section of $p_\#$.  It is clear
  that $s_\#$ is preserved under the bubblings of $s$.
  It is easy to check that   the resulting mapping
    \begin{equation}\label{9.2.a}{\mathcal S}_*(p)\to S_*(p_\#)\, ,\,\,\,  s\mapsto s_\#\end{equation} is a bijection.
   Thus, $\vert {\mathcal S}_*(p)\vert=\vert {  S}_*(p)\vert$.     Corollary \ref{1.3.4} implies the
    following claim.

   \begin{theor}\label{9.2.1}   If
  $\Phi $ is finite, then
\begin{equation}\label{9.2.b}\vert {\mathcal S}_*(p)\vert=\vert  \Phi \vert \,
\sum_{\rho\in \Irr(\Phi),\, \pi_\rho = \pi   } (\vert  \Phi
\vert/\dim \, \rho)^{-\chi(W)}\, \zeta_\rho  ([W]) \, .
\end{equation}
\end{theor}

  We now rewrite (\ref{9.2.b}) in terms of non-pointed sections
of $p$. Since $F$ is path-connected,  any section of $p$ is
homotopic to a pointed section. This shows that the natural mapping
${\mathcal S}_*(p)\to {\mathcal S}(p)$ is surjective. This mapping
may be described in terms of an action of $\Phi$ on ${\mathcal
S}_*(p)$ as follows. The group $\Phi$ acts on the set $S_*(p_\#)$ by
conjugation. This defines an action of $\Phi$ on ${\mathcal S}_*(p)$
via the bijection (\ref{9.2.a}). It is easy to see that the orbits
of the latter action are precisely the preimages of elements of
${\mathcal S}(p)$ under the natural mapping ${\mathcal S}_*(p)\to
{\mathcal S}(p)$. Thus, ${\mathcal S}(p)={\mathcal S}_*(p)/\Phi$.
For $s\in {\mathcal S}(p)$, let $\Aut(s)\subset \Phi$ be the
stabilizer
  of an element of ${\mathcal S}_*(p)$ projecting to $s$. The group $\Aut(s)$
is well defined up to conjugation in  $\Phi$. If $\Phi$ is finite,
then
$$\vert {\mathcal S}_*(p)\vert=  \sum_{s\in {\mathcal S}(p)} \, \frac{\vert \Phi\vert} {\vert \Aut(s)
 \vert}\, .$$
 Formula (\ref{9.2.b}) may now be rewritten as
\begin{equation}\label{9.2.c}\sum_{s\in {\mathcal S}(p)} \, \frac{1} {\vert \Aut(s)
 \vert}=\sum_{\rho\in \Irr(\Phi),\, \pi_\rho = \pi   } (\vert  \Phi \vert/\dim \,
\rho)^{-\chi(W)}\, \zeta_\rho  ([W])\, . \end{equation}

   Theorem  \ref{1.3.2} and Corollary \ref{1.3.3} imply that
    $\sum_{s\in {\mathcal S}(p)} \,  {1} /{\vert \Aut(s)
 \vert}$ is an integer
     divisible by $\vert Z(\Phi)\vert^{2d-2}$,
   where   $d$ is the genus of~$W$,
   and
   $$\sum_{s\in {\mathcal S}(p)} \,  {1} /{\vert \Aut(s)
 \vert}\,\, \leq \,\,\vert \Phi\vert^{-1}\,
   \vert \Hom (\pi, \Phi)\vert\, .$$

 \begin{corol}\label{9.2.1+} Under the assumptions of Theorem \ref{9.2.1}
the bundle $p:E\to W$ has a section if and only if $$\sum_{\rho\in
\Irr(\Phi),\, \pi_\rho = \pi   } ( \dim \, \rho)^{\chi(W)}\,
\zeta_\rho  ([W]) \neq 0\, .$$
\end{corol}

Note that the left-hand side of this formula is a non-negative
rational number for all $p$, as it follows from  Theorem
\ref{9.2.1}. In the case where $\Phi$ is abelian or the genus of $W$
is bigger than $(1/2)\,\log_2(\vert
 [\Phi, \pi']\vert -1)$, Theorems \ref{7.1.1} and \ref{7.1.2} imply that the bundle $p$
  has a section if and only if the induced
homomorphism $p_*:H_2(E)\to H_2(W)$ is surjective.

In the case of a trivial fiber bundle, Theorem \ref{9.2.1} amounts
to computing the number of pointed homotopy classes of maps $W\to
F$. In this case, all the cohomology classes $\zeta_\rho$ are
trivial and Theorem \ref{9.2.1} follows   from the equality $\vert
{\mathcal S}_*(p)\vert=\vert {  S}_*(p)\vert$ and  the
Frobenius-Mednykh formula (\ref{1.1.a}).

  \subsection* {Enumeration of lifts of maps.} Given an arbitrary locally
trivial fiber bundle $p:E\to X$ and a map $f$ from the surface $W$
to   $X$, one may be interested in counting  the number of homotopy
classes of lifts of $f$ to $E$.   By a {\it lift}\index{lift} of $f$
to $E$, we mean a mapping $f':W\to E$ such that   $pf'=f$. For $X=W$
and $f=\id_W$, we recover the setting of the previous subsection.
All definitions and results given there extend to arbitrary $p,f$
with the obvious changes. The key observation is that the lifts of
$f$ to $E$ bijectively correspond to the sections of the induced
fiber bundle $f^*(p)$ over $W$.

\subsection* {Non-abelian cohomology of  surfaces.}   Theorem \ref{t1} yields interesting information about
1-dimensional non-abelian cohomology of the fundamental groups of
surfaces. We begin by
  recalling the definition of the 1-dimensional non-abelian
cohomology of an arbitrary group $\pi$, cf.\ [Se]. Fix a left action
of $\pi$
 on a group $\Phi$, i.e., a homomorphism $\pi\to \Aut \, \Phi$. A map $\alpha:\pi\to \Phi$
is  a {\it cocycle}\index{cocycle} if $\alpha (ab)=\alpha (a)\,
a(\alpha(b))$ for all $a,b \in \pi$. Here $a(\alpha(b))\in \Phi$ is
obtained by  the action of $a$ on $\alpha(b)$. For example, the
mapping $\pi\to \{1\}\subset \Phi$ is a cocycle. The set of all
cocycles $\pi\to \Phi$  is denoted by $Z^1(\pi;\Phi)$. The group
$\Phi$ acts on $Z^1(\pi;\Phi)$ by
$$(\varphi \alpha) (a)= \varphi \, \alpha(a) \, (a\varphi)^{-1}$$
for all $\varphi\in \Phi$, $\alpha\in Z^1(\pi;\Phi)$, and $a\in
\pi$. The quotient set of this action is denoted by $H^1(\pi;\Phi)$
and called the (nonabelian) cohomology of $\pi$ with coefficients in
$\Phi$.  For   $h\in H^1(\pi;\Phi)$, let $\Aut (h)\subset \Phi$ be
the stabilizer   of any cocycle representing $h$. The group $\Aut
(h) $ is well defined  up to conjugation in $\Phi$.

If   $\pi$ is finitely generated and $\Phi$ is finite, then both
sets $Z^1(\pi;\Phi)$ and $H^1(\pi;\Phi)$ are finite. Put
$$\mathcal M  (\pi;\Phi)\, =\,  \sum_{h\in H^1(\pi;\Phi)} \,\frac{1 }{\vert \Aut (h )
\vert}\in \QQ \, .$$ We view $\mathcal M  (\pi;\Phi)$ as the global
measure of the set $H^1(\pi;\Phi)$ counting its   elements   with
the weights
  $ {1 }/{\vert \Aut   \vert}$. Since   any $h\in
H^1(\pi;\Phi)$ can be
  represented by precisely $\vert \Phi\vert/\vert \Aut (h ) \vert$
  cocycles,
  $$\mathcal M  (\pi;\Phi)\, =\,\vert \Phi\vert^{-1}\, \vert Z^1(\pi;\Phi)\vert\, .$$

The definitions of Section  \ref{eh}  can be adapted to this setting
as follows. With an   irreducible  representation $\rho:\Phi\to
GL_n(\CC)$ of $\Phi$ we associate the group $\pi_\rho\subset \pi$
consisting of all $a\in \pi$ such that the representation
$\varphi\mapsto \rho (a^{-1}\varphi)$ of $\Phi$ is equivalent
to~$\rho$. This means that there is a matrix $M_a\in GL_n(\CC)$ such
that $\rho (a^{-1}\varphi)=M_a^{-1}\,  \rho(\varphi) \, M_a$ for all
$\varphi\in \Phi$. Then there is a family of non-zero complex
numbers $\{\zeta_{a,b} \}_{a,b\in \pi_\rho}$ such that
$\zeta_{a,b}\, M_a \, M_b=M_{ab}$ for all $a,b\in \pi_\rho$. This
family is
  a 2-cocycle representing a well-defined cohomology class
$\zeta_\rho\in H^2(\pi_\rho; \CC^*)$.

  \begin{theor}\label{8.4.1}     For any action
of $\pi=\pi_1(W)$ on a finite group $\Phi$,
$$\mathcal M  (\pi;\Phi)\,
  =\,
\sum_{\rho\in \Irr(\Phi),\, \pi_\rho =\pi   } (\vert  \Phi
\vert/\dim \, \rho)^{-\chi(W)}\,\, \zeta_\rho  ([W]) \,.  $$
\end{theor}

    \begin{proof}
  Let $\pi'$ be the   set of  pairs $(\varphi\in \Phi, a\in \pi)$
   with multiplication $(\varphi , a)(\varphi' , a')=(\varphi \, (a\varphi') , aa')$.
   It is easy to check that $\pi'$ is a group. The formula $p(\varphi ,
   a)= a$ defines an epimorphism $p:\pi'\to \pi$ with kernel
   $\{(\varphi, 1)\}_{\varphi\in \Phi}=\Phi$.
   Every cocycle $\alpha:\pi\to \Phi$ defines a section $s_\alpha$ of $p$  by
   $s_\alpha(a)= (\alpha (a), a)$ for $a\in \pi$. The formula $\alpha\mapsto s_\alpha$  establishes a
   bijection between the set   $Z^1(\pi;\Phi)$ and the set
    $S_*(p)$ of the sections of $p$. Therefore
    \begin{equation}\label{8.4.aa}\mathcal M  (\pi;\Phi)\, =
    \,\vert \Phi\vert^{-1}\, \vert Z^1(\pi;\Phi)\vert\, =\, \vert \Phi\vert^{-1}\,  \vert
     S_*(p)\vert\,.\end{equation}
    It remains to apply Corollary \ref{1.3.4} and to
   observe that  the definitions of $\pi_\rho$ and $\zeta_\rho$ given
   in Section \ref{eh} are equivalent in the present setting to the definitions
    given before the statement of the theorem. (The key point is that every $a\in \pi$ has a canonical lift $(1,a)$ to $\pi'$ and
     $$(1,a)^{-1} (\varphi, 1)(1,a)=(a^{-1} \varphi, 1)$$ for all   $\varphi\in \Phi$.)
     \end{proof}

     Formula (\ref{8.4.aa}) and the remarks   after Corollary \ref{1.3.4} imply that
    $\mathcal M  (\pi;\Phi)$ is an integer divisible by $\vert Z(\Phi)\vert^{2d-2}$,
   where    $d$ is the genus of $W$, and
   $$\mathcal M  (\pi;\Phi) \leq \vert \Phi\vert^{-1}\,
   \vert\Hom (\pi, \Phi)\vert\, .$$

\section {Miscellaneous algebra}

We discuss   miscellaneous  algebraic notions and results related to
  Theorem \ref{t1}.

 \subsection* {Extremal homology classes.}
 We call a homology class $h\in H_2(G{{}})$ {\it
 extremal}\index{homology class!extremal} (with respect to the given epimorphism $q:G'\to G$)  if
  $ \zeta_\rho  (h)=1$ for all irreducible
    representations $\rho$ of $\Gamma=\Ker\, q$  such that $G_\rho=G$.
    For example, the zero homology class $h=0$ is extremal. It is
    clear that the extremal homology classes form a
    subgroup of~$H_2(G{{}})$.

    For each $k\geq 1$, the function  $v_k$ introduced in Section \ref{vvv}    takes
    on all
    extremal classes the  same value which is the maximal value of
    $v_k$.
    In particular, if $h \in H_2(G{{}})$ is extremal, then
    $v_1(h)=v_1(0)>0$. By Lemma \ref{l5}(a),   all extremal homology classes lie in
  $Q=\Im\, (q_*: H_2(G'{{}})\to H_2(G{{}}))$. Let $a\in \ZZ$ be the least common multiple of the
numbers $ \dim \, \rho$, where $\rho$ runs over all irreducible
  representations of~$\Gamma$.   The properties of $\zeta_\rho$ imply that  all
elements of the group $a   Q\subset Q$ are extremal.

\subsection* {Quasi-epimorphisms.}   A homomorphism
$g:\pi \to G$ is a  {\it
quasi-epimorphism}\index{quasi-epimorphisms} (with respect to
$q:G'\to G$) if $$g(\pi)\cap (G-G_\rho)\neq \emptyset$$ for all
irreducible   representations $\rho$ of $\Gamma$ such that
$G_\rho\neq G$. In particular, all epimorphisms $\pi\to G$ are
quasi-epimorphisms. If $G_\rho= G$ for all $\rho$, then all
homomorphisms $\pi\to G$ are quasi-epimorphisms.

  For a    quasi-epimorphism
$g:\pi=\pi_1(W)\to G$, Theorem \ref{t1} yields the same formula
(\ref{1.3.a+}) as for an epimorphism.
  Therefore
$$\vert  \Hom_g (\pi ,G')\vert\,\leq\, \vert  \Gamma \vert  \,
\sum_{\rho\in \Irr(\Gamma),\, G_\rho =G  } (\vert  \Gamma \vert/\dim
\, \rho)^{-\chi(W)}\,  .  $$ This  inequality is an equality if and
only if  $g_*([W])\in H_2(G{{}})$ is an extremal homology class.

\subsection* {The genus norm.} As was already mentioned above, for
any non-zero $h\in H_2(G{{}})$, there are a closed connected
oriented surface $\Sigma$ of positive genus and a homomorphism
$g:\pi_1(\Sigma)\to  G$ such that
  $g_*([\Sigma])=h$.
Let $\vert h\vert\geq 1 $   be the minimal
  genus of such   $\Sigma$. By definition, $\vert 0\vert =0$.
   Clearly,   $\vert h+h'\vert \leq \vert h\vert+ \vert h'\vert$ for
all $h, h'\in H_2(G{{}})$. Also, $\vert -h\vert=\vert h\vert  $ and
$\vert h\vert =0$ if and only if $h=0$.
  We call the   mapping $H_2(G{{}})\to \ZZ,\, h\mapsto \vert h\vert$
  the {\it genus norm}\index{genus norm}.

  The functions $v_1, v_2,
\ldots$   from  Section \ref{vvv} may help to
  estimate the genus norm as follows.
Define a mapping $v: H_2(G{{}})\times \ZZ \to \QQ$ by
$$v(h,n)=\sum_{k\geq 1  }  \,v_k (h)\,k^{ - 2n}\, , $$
where  $h\in H_2(G{{}})$ and $n\in \ZZ$.

  \begin{lemma}\label{l7}  If $G_\rho= G$
for all irreducible   representations $\rho$ of $\Gamma$, then
$v(h,n)\geq 0$ for all $h\in H_2 (G{{}})$ and all  $n\geq  \vert
h\vert-1$.\end{lemma}

    \begin{proof} If $h=0$,
 then $v(h,n)> 0$ for all   $n\in \ZZ$, as directly follows from
  Lemma \ref{l5}(c). Suppose that $h\neq 0$. Let $\Sigma$ be  a closed connected oriented  surface  of genus
 $\vert h\vert\geq 1  $
  and let  $g:\pi_1(\Sigma) \to  G$ be  a
homomorphism  such that
  $g_*([\Sigma])=h$. Theorem \ref{t1} and the assumption $G_\rho= G$
for all $\rho$ imply that
$$\vert\Hom_g (\pi_1(\Sigma) ,G') \vert    =   \vert  \Gamma \vert^{1-\chi(\Sigma)}
 \, v (h, -\frac{\chi(\Sigma)}{2})=\vert  \Gamma \vert^{ 2\vert h\vert-1}
 \, v (h,  \vert h\vert-1)\,  .  $$
Hence, $v(h, \vert h\vert-1)\geq 0$.  Similarly $v(h,n)\geq 0$ for
all   $n\geq  \vert h\vert$,  cf.\ the argument  after
  (\ref{vv2}).    \end{proof}

  For  $h\in H_2 (G{{}})$,   denote by   $\langle h,
q\rangle  $  the minimal non-negative integer such that $v(h,n)\geq
0$ for all   $n\geq  \langle h, q\rangle$. Lemma \ref{l7} implies
that such an integer  exists and   $\vert h\vert
 \geq \langle h,
q\rangle +1$. Varying $q $ in the class of group epimorphisms
 with target~$G$ and finite kernel satisfying the conditions of Lemma \ref{l7}, we obtain a
  family of estimates  from
  below for the genus norm on~$H_2(G{{}})$.
 The author does not know whether these estimates may be
 non-trivial. Explicit examples and    computations would be welcome.

 \section {A generalization   of   Theorem \ref{t1}}\label{sect2}

Theorem \ref{t1} admits a  generalization  in which the lifts of $g$
to  $ G'$ are counted with weights determined by a 2-cocycle on
$G'$. To state this generalization, we first extend the definition
of $\zeta_\rho$ to projective representations of $\Gamma$.

Fix  throughout this section  a   2-cocycle $\theta=\{\theta_{a,b}
\in \CC^*\}_{a,b\in \orlin'}$ on $G'$.  A mapping $\rho:\Gamma \to
GL_n (\CC)$ with $n=1,2,\ldots$ is   a {\it $\theta$-representation
of}\index{$\theta$-representation} $\Gamma$
  if $\rho(1) $ is the unit $n\times n$ matrix  and $  \rho (g) \,\rho (h)\, =\, \theta_{g, h}\, \rho(gh) $
   for all $g,h\in \Gamma$. For any  matrix $M\in GL_n(\CC)$ and a $\theta$-representation $\rho:\Gamma \to GL_n (\CC)$,
     the mapping $M^{-1}\rho \,M:\Gamma \to GL_n (\CC)$ sending   $h\in \Gamma$ to $M^{-1}  \rho (h)  M$, is
    a $\theta$-representation.
     We say that two $\theta$-representations $\rho :\Gamma\to GL_{n} (\CC)$ and
 $\rho':\Gamma\to GL_{n'}(\CC)$ are {\it equivalent}
 and write $\rho \sim \rho'$ if $ n =n'$ and   $\rho' =M^{-1}\rho  \,  M$ for some $M\in GL_{n }(\CC)$.  Clearly, $\sim$ is an equivalence relation on the set of
 $\theta$-representations of $\Gamma$.
Denote by ${\mathcal R}_\theta$ the corresponding set of equivalence
classes.

 The cocycle  $\theta$ determines an
 action of $G'$ on ${\mathcal R}_\theta$ as follows.
 Given     $a\in G'$ and a $\theta$-representation $\rho:\Gamma\to GL_n(\CC)$,
 consider the mapping ${ a\rho}: \Gamma\to
  GL_n(\CC)$, whose value on any $h\in \Gamma$ is given by
    $${ a\rho}(h)=  \frac { \theta_{ a^{-1},  h
    a }\, \theta_{   h, a }  }{ \theta_{ a,  a^{-1} } \,  \theta_{ 1,  1 }  \,} \, \rho ( a^{-1} h
  a ) \,  .  $$

   \begin{lemma}\label{3.1.1} The mapping ${ a\rho}$ is a $\theta$-representation of $\Gamma$.
   The formula $(a, \rho)\mapsto {a\rho}$ defines a  left  action of $G'$
   on~${\mathcal R}_\theta$. This action induces a left action of    $G$ on~${\mathcal R}_\theta$.
\end{lemma}

Given a $\theta$-representation $\rho$ of $\Gamma$, denote by
  $G_\rho $ be the stabilizer of $\rho$, i.e.,  the subgroup of $G$ consisting
  of
  all $\alpha\in   \orlin$ such that  $a\rho\sim \rho$ for
 some (and then for all) $a\in q^{-1}(\alpha)\subset G'$. The group $G_\rho$ depends  only on
   the equivalence class of~$\rho$.

Let $\rho:\Gamma\to GL_n(\CC)$ be a   $\theta$-representation of
$\Gamma$,
  which is irreducible in the sense
 that   the only linear subspaces of $\CC^n$ preserved
under the induced projective action of $\Gamma$ are $0$ and $\CC^n$.
 We   define a cohomology class $\zeta_\rho\in H^2(G_\rho; \CC^*)$. Fix for each $\alpha\in
G_\rho $, an element $\widetilde \alpha$ of $
 q^{-1} (\alpha)$.
  By definition of
 $G_\rho$, for $\alpha\in G_\rho$,     there is a matrix $M_\alpha\in GL_n(\CC)$ such that
 $  {\widetilde \alpha}   \rho   = M_\alpha^{-1} \, \rho   \,
 M_\alpha$.
 The irreducibility of $\rho$ implies that  $M_\alpha$ is
  unique  up to multiplication by an element of $\CC^*$; we fix $M_\alpha$ for all $\alpha$.

\begin{lemma}\label{3.3.1}    For any  $\alpha,
\beta\in G_\rho$, there is a unique   $\zeta_{\alpha, \beta} \in
\CC^*$ such that $$\zeta_{\alpha, \beta} \,M_\alpha\,  M_\beta=
 \theta_{\widetilde \alpha, \widetilde \beta}\,
\theta^{-1}_{ \widetilde {\alpha\beta}, \widetilde
{\alpha\beta}^{-1}\, \widetilde \alpha \, \widetilde \beta}\,
M_{\alpha \beta} \, \rho ( \widetilde {\alpha\beta}^{-1}\,
\widetilde \alpha \, \widetilde \beta) .  $$ The family
$\{\zeta_{\alpha,\beta}\}_{\alpha,\beta}$ is a   2-cocycle on
 $G_\rho$. Its cohomology class  $\zeta_\rho\in H^2(G_\rho; \CC^*)$    depends only
  on the equivalence
class of $\rho$ and does not depend on the choice
 of the matrices $\{M_\alpha\}_{\alpha}$  or the lifts $\{\widetilde \alpha\}_{\alpha }$.
 \end{lemma}

There is a relationship between   $\zeta_\rho$ and the cohomology
class of the cocycle $\theta$. Consider the group
$G'_\rho=q^{-1}(G_\rho)\subset G'$ and let
 $[\theta]_\rho\in H^2( G'_\rho;\CC^*)$
denote  the cohomology class of the restriction of   $\theta$ to
$G'_\rho$. In generalization of (\ref{1.3.ddd}), we have
\begin{equation}\label{nn}(\dim\, \rho) \, q^*(\zeta_\rho)=(\dim\,
\rho)\, [\theta]_\rho\, \end{equation}  where
$q^*:H^2(G_\rho;\CC^*)\to H^2( G'_\rho;\CC^*)$ is the homomorphism
induced by~$q$.

  \begin{theor}\label{t2}     Let
$\theta=\{\theta_{a,b} \in \CC^*\}_{a,b\in G'}$ be a  2-cocycle
  on~$G'$ representing a cohomology class $[\theta]\in
H^2(G'; \CC^*)$. For any group homomorphism $g:\pi=\pi_1(W)\to G$,
\begin{equation*}  \sum_{g'\in \Hom_g (\pi ,G')} (g')^* ([\theta])([W])
\, =\, \vert  \Gamma \vert  \, \sum_{\rho,\, G_\rho \supset \,g(\pi)
} (\vert  \Gamma \vert/\dim \, \rho)^{-\chi(W)}\,(g^*(\zeta_\rho)
([W]))^{-1} \, ,
\end{equation*} where $\rho$ ranges over the equivalence classes of
irreducible
  $\theta$-representations of $\Gamma$
such that $G_\rho \supset g(\pi)$.
 \end{theor}

     Here   $(g')^*
([\theta]) ([W])\in \CC^*$ and  $g^*(\zeta_\rho) ([W]) \in \CC^*$
are the evaluations of   $(g')^* ([\theta])$, $ g^*(\zeta_\rho) \in
H^2(\pi; \CC^*) $ on $[W]\in H_2(W{{}})=H_2(\pi{{}})$, respectively.

Applying  Theorem \ref{t2}    to $-W$, we obtain
$$ \sum_{g'\in \Hom_g (\pi ,G')}
(g')^* ([\theta])([-W]) =\, \vert  \Gamma \vert  \, \sum_{\rho,\,
G_\rho \supset \,g(\pi)  } (\vert  \Gamma \vert/\dim \,
\rho)^{-\chi(W)}\,g^*(\zeta_\rho) ([W]) \, .  $$ For $\theta=1$, the
left-hand side   is equal to $\vert \Hom_g (\pi ,G')\vert$, and we
obtain   Theorem \ref{t1}. For  arbitrary $\theta$ and $g=1$,
Theorem \ref{t2} is   contained in \cite{Tu}. A proof of Theorem
\ref{t2} in the general case will be given elsewhere.

Theorem \ref{t2} yields a generalization of other formulas obtained
above. We state   generalizations of Formulas (\ref{ffmc}) and
(\ref{9.2.c}) using notation of Section \ref{kol}.
 For
a   principal $G'$-bundle $\xi'$ on $W$, set $$ \theta_{\xi'}=
(g')^*([\theta]) \in H^2(\pi; \CC^*),\, $$ where $[\theta]\in
H^2(G'; \CC^*)$ is the   cohomology class of $\theta$   and
$g':\pi\to G'$ is any homomorphism such that $\xi'=\xi_{g'}$. The
cohomology class $\theta_{\xi'}$ does not depend on the choice of
$g'$ because the conjugations in $G'$ act trivially on
$H^*(G';\CC^*)$.   Theorem \ref{t2}  implies that for the
principal $G$-bundle $\xi$ over~$W$ determined by a homomorphism
$g:\pi\to G$,
$$\sum_{\xi'\in {\mathcal P}(\xi)} \,\frac{\theta_{\xi'}([W])}{\vert \Aut_\xi (\xi' ) \vert} \,
  =\,
\sum_{\rho,\, G_\rho \supset \,g(\pi)  } \left ( \frac {\vert \Gamma
\vert} {\dim \, \rho}\right )^{-\chi(W)}\,\,(g^*(\zeta_\rho)
([W]))^{-1}\, ,
$$ where $\rho$ ranges over the equivalence classes of irreducible
   $\theta$-representations of $\Gamma$ such that
$G_\rho \supset g(\pi)$. For $G=1$,  this formula boils down to
$$\sum_{\xi'\in {\mathcal P} } \,\frac{\theta_{\xi'}([W])}{\vert \Aut  (\xi' ) \vert} \,
  =\,
\sum_{\rho  } \left ( \frac {\vert \Gamma \vert} {\dim \,
\rho}\right )^{-\chi(W)}\,  ,
$$ where $\rho$ ranges over the equivalence classes of irreducible
   $\theta$-representations of $\Gamma$.

  To   generalize   (\ref{9.2.c}), consider a locally trivial fiber
  bundle $p:E\to W$ and
  a cohomology class $\Theta\in H^2(E;\CC^*)$ whose
evaluation on $\pi_2(E)$ is equal to $1 $. Such $\Theta$ is
necessarily induced from a unique   element of $
H^2(\pi_1(E);\CC^*)$. We represent the latter by a $\CC^*$-valued
2-cocycle $\theta$ on $ \pi_1(E )$. Set $\pi=\pi_1(W)$.

   \begin{theor}\label{9.2.2} If the fiber of $p$ is path-connected
   and has a finite fundamental group
  $\Phi $, then
$$\sum_{s\in {\mathcal S}(p)} \, \frac{s^*(\Theta)([W])} {\vert \Aut(s)
 \vert}=\sum_{\rho,\, \pi_\rho = \pi   } (\vert  \Phi \vert/\dim \,
\rho)^{-\chi(W)}\, (\zeta_\rho  ([W]))^{-1}\, ,  $$ where $\rho$
runs over the equivalence classes of irreducible
   $\theta$-representations of~$\Phi$   such that
$\pi_\rho =\pi $.
\end{theor}

  This   follows from Theorem \ref{t2}  using the
arguments   of Section \ref{kol}. The assumption
$\Theta(\pi_2(E))=1$
  ensures that $s^*(\Theta)\in H^2(W;\CC^*)$ is a bubble
equivalence invariant of  a section $s$ so that its evaluation on $
[W]$ is well-defined.

                     \end{document}